\documentclass[reqno,11pt]{amsart}

\usepackage{amsmath,amssymb}
\usepackage{enumerate}
\usepackage{graphicx}
\usepackage[hidelinks]{hyperref}

\theoremstyle{plain}
\newtheorem{thm}[equation]{Theorem}
\newtheorem{cor}[equation]{Corollary}
\newtheorem{lem}[equation]{Lemma}
\newtheorem{prop}[equation]{Proposition}

\theoremstyle{remark}
\newtheorem{rmk}[equation]{Remark}

\numberwithin{equation}{section}

\newcommand{\CC}{\mathbb{C}}
\newcommand{\HH}{\mathbb{H}}
\newcommand{\RR}{\mathbb{R}}
\newcommand{\scrE}{\mathcal{E}}
\newcommand{\scrL}{\mathcal{L}}
\newcommand{\scrM}{\mathcal{M}}
\newcommand{\scrP}{\mathcal{P}}
\newcommand{\scrS}{\mathcal{S}}
\newcommand{\deGa}{\delta_{\Gamma}}

\DeclareMathOperator{\PSL}{PSL}
\DeclareMathOperator{\Ree}{Re}
\DeclareMathOperator{\Res}{Res}
\DeclareMathOperator{\dist}{dist}
\DeclareMathOperator{\htop}{ht}
\DeclareMathOperator{\spec}{spec}
\DeclareMathOperator{\supp}{supp}
\DeclareMathOperator{\vol}{vol}
\newcommand{\Stab}{\mathrm{Stab}}
\newcommand{\euc}{\mathrm{Eucl}}
\newcommand{\norm}[1]{\lVert#1\rVert}
\newcommand{\dd}{\,d}

\title[Counting Circles]{Counting circles in Kleinian circle packings}
\author{Samuel C. Edwards}
\address{Department of Mathematics, Uppsala University,
  Box 524, 751 20 Uppsala, Sweden}
\email{samuel.edwards@math.uu.se}
\date{}

\begin{document}

\begin{abstract}
Let $\mathcal{P}$ be a bounded Kleinian circle packing with symmetry group
$\Gamma$, and let $N_{\mathcal{P}}(T)$ denote the number of circles in
$\mathcal{P}$ of radius at least $T^{-1}$. We obtain an asymptotic formula for
$N_{\mathcal{P}}(T)$ in which the eigenvalues of the Laplace operator below the continuous
spectrum contribute a finite sum of main terms, with an error term having the same exponent as
in the counting theorem of Lax--Phillips for
$\Gamma$-orbits in $\mathbb H^3$.

In the proof, we express the Dirichlet series associated to $N_{\mathcal{P}}(T)$
as an integral of an Eisenstein series over a totally geodesic surface in
$\Gamma\backslash\mathbb{H}^3$. We then use the spectral theorem to
meromorphically continue this Dirichlet series, and a point-pair kernel to
obtain the spectral bounds needed for counting.
\end{abstract}
\maketitle

\section{Introduction}
Let $\scrP$ be a bounded Kleinian circle packing (cf.~\cite{KK,KMS,Parker}).
This means that $\scrP$ is an infinite collection of circles $C$ in
$\widehat{\CC}=\CC\cup\{\infty\}$ with a choice of round open disk $D_C$ bounded
by each $C$, such that
\begin{enumerate}[(1)]
\item the disks $D_C$, $C\in\scrP$, are pairwise disjoint and their union is
  dense in $\widehat{\CC}$;
\item the residual set $\Lambda_{\scrP}:=\widehat{\CC}-\bigcup_{C\in\scrP}D_C$
  is the limit set $\Lambda_\Gamma$ of a geometrically finite Kleinian group
  $\Gamma<\PSL(2,\CC)$, which we call a symmetry group of $\scrP$;
\item $\infty$ is inside one of the disks $D_C$.
\end{enumerate}

\begin{figure}[b]
\begin{center}
\includegraphics[width=.70\textwidth]{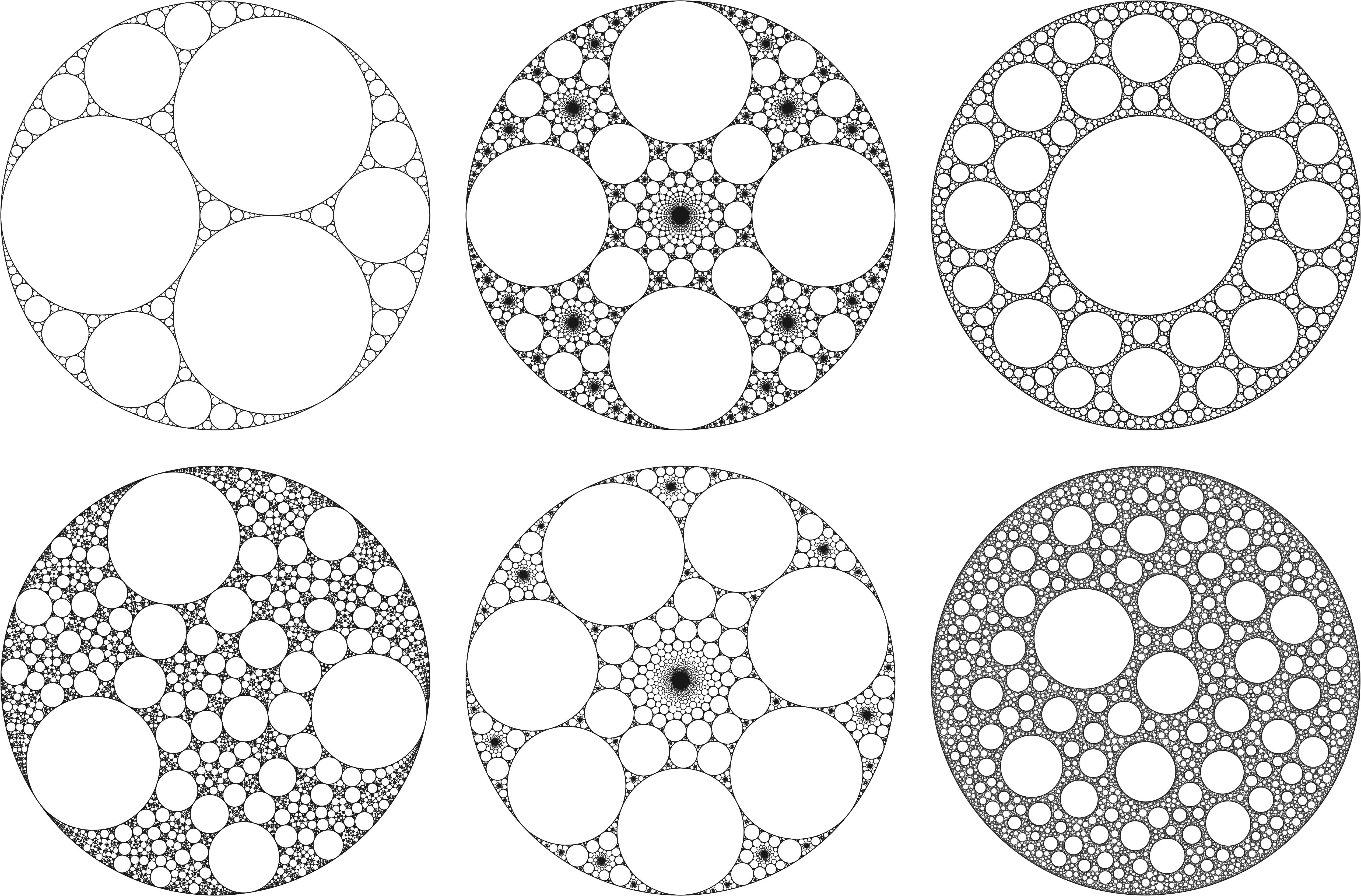}
\caption{Examples of bounded Kleinian circle packings.}
\label{fig:packings}
\end{center}
\end{figure}

For a circle $C\in\scrP$, let $\mathsf{r}(C)>0$ denote its Euclidean radius. We
are interested in the counting function
\begin{equation}
  \label{eq:counting}
  N_{\scrP}(T):=\#\{C\in\scrP:\mathsf{r}(C)\ge T^{-1}\},\qquad  T>0.
\end{equation}

We study $N_{\scrP}(T)$ as $T\to\infty$ using the Dirichlet series
\begin{equation}
  \label{eq:zeta}
  Z_{\scrP}(s)=\sum_{C\in\scrP}\mathsf{r}(C)^s=s\int_0^\infty N_{\scrP}(T)T^{-s-1}\dd T.
\end{equation}
The disks not containing $\infty$ are disjoint and lie in a bounded region,
so $\sum_{C\in\scrP}\mathsf r(C)^2<\infty$.
The series and the integral are therefore absolutely convergent for $\Ree s\geq2$. We will meromorphically continue $Z_{\scrP}(s)$ to the half-plane
$\Ree s>1$ and bound it near the line $\Ree s=1$.

For each circle $C\in \scrP$, let $\Pi_C\subset\HH^3$ denote the hyperbolic
plane whose ideal boundary is $C$. The disks $D_C$ are the components of the
domain of discontinuity $\Omega(\Gamma)$ of $\Gamma$. Ahlfors' finiteness
theorem \cite{Ahlfors} implies that $\Omega(\Gamma)$ consists of finitely many
$\Gamma$-orbits of these components, and that each quotient
$\Gamma_C\backslash D_C$, where $\Gamma_C=\Stab_\Gamma(C)$, has finite
hyperbolic area. Identifying $D_C$ with $\Pi_C$, we choose representatives
$C_1,\ldots,C_n$ of these orbits, so that $\scrP=\bigsqcup_{i=1}^n\Gamma C_i$,
and define
\[
  \scrS=\bigsqcup_{i=1}^n\Gamma_{C_i}\backslash\Pi_{C_i}.
\]
Then $\scrS$ is a finite-area totally geodesic surface in
$\scrM=\Gamma\backslash\HH^3$, possibly with several components and orbifold
points. We write $dS$ for the hyperbolic area measure on the planes and on
$\scrS$. The quotient $\scrM$ is a geometrically finite hyperbolic $3$-orbifold
of infinite volume.

We introduce the following notation for the upper half-space model of $\HH^3$:
\[
  p=z+rj,\qquad z\in\CC,\ r>0,
  \qquad \htop(p)=r,
  \qquad \dd\vol(p)=\frac{dx\,dy\,dr}{r^3}.
\]
Let $\deGa$ be the critical exponent of $\Gamma$. For any $C\in\scrP$,
the stabilizer $\Gamma_C$ is a lattice in the isometry group of $\Pi_C$. In
particular, $\delta_{\Gamma_C}=1$ and $\Gamma_C$ is of divergence type. Since
$\Lambda_{\Gamma_C}=C\subsetneq\Lambda_\Gamma$, Dal'bo--Otal--Peign\'e
\cite[Proposition~2]{DOP} gives $\deGa>\delta_{\Gamma_C}=1$. On
the other hand, $\deGa<2$ since $\Gamma$ is geometrically finite and has
infinite covolume (see \cite[Theorem~1]{Sullivan}). We define an Eisenstein
series $E(\cdot,s)$ on $\HH^3$ by
\begin{equation}
  \label{eq:eisenstein}
  E(p,s)=\sum_{\gamma\in\Gamma}\htop(\gamma p)^s,\qquad p\in\HH^3,\; \Ree s>\deGa .
\end{equation}
The series is absolutely convergent for $\Ree s>\deGa$ (cf., e.g.,
\cite[Proof of Theorem 5.2]{EO}).

Letting $\scrL$ denote the negative of the Laplace operator on $\HH^3$ or
$\scrM$, we have
\[
  \scrL E(p,s)= \lambda_sE(p,s),\qquad \text{where }\lambda_s=s(2-s).
\]

For $\Ree s>1$, integration over $\Pi_C$ gives
\[
  \int_{\Pi_C}\htop(p)^s\dd S(p)
  =\frac{2\pi}{s-1}\mathsf{r}(C)^s.
\]
Summing over the packing planes for $\Ree s\geq2$ gives
\begin{equation}
  \label{eq:period}
  Z_{\scrP}(s)=\frac{s-1}{2\pi}
  \int_{\scrS}E(p,s)\dd S(p).
\end{equation}
Bunke and Olbrich \cite{BunkeOlbrich} and Guillarmou and Mazzeo
\cite{GuillarmouMazzeo} proved meromorphic continuation of $E(p,s)$
(see also \cite{FHP,MazzeoMelrose,Perry}). Since $\scrS$ need not be compact, their results do not
directly give meromorphic continuation or bounds for the integral in
\eqref{eq:period}. We use the spectral theorem to meromorphically continue
$Z_{\scrP}(s)$, and the spectral measure of the surface measure on $\scrS$ to
bound it on vertical lines.

By Sullivan \cite{SullivanPositivity}, the base
eigenvalue of $\scrL$ on $L^2(\scrM)$ is $\deGa(2-\deGa)$. It is simple, with
positive eigenfunction $\phi_0\in L^2(\scrM)$. Lax and Phillips
\cite{LaxPhillips} proved that there are only finitely many eigenvalues below
the continuous spectrum $[1,\infty)$. We write the remaining eigenvalues below
$1$ as $s_j(2-s_j)$, where
\[
  2>\deGa>s_1\ge s_2\ge\cdots\ge s_J>1.
\]
If there is no eigenvalue in $(\deGa(2-\deGa),1)$, we let $J=0$ and the list
$s_1,\ldots,s_J$ is empty.

\begin{thm}
\label{thm:main}
\leavevmode
\begin{enumerate}[(i)]
\item Let $\scrE=\{\deGa\}\cup\{s_j:1\le j\le J\}$. The function $Z_{\scrP}(s)$
  has a meromorphic continuation to the half-plane $\Ree s>1$. Its only possible
  poles are at the points of $\scrE$, and are at most simple. The residue at
  $s=\deGa$ is positive.

\item For $1<\sigma\le2$ and $|t|\ge1$,
\[
  |Z_{\scrP}(\sigma+it)|
  \ll_{\scrP,\Gamma}
  \frac{|t|}{\sigma-1}.
\]
Moreover, there exists $\varepsilon_0>0$ such that
\[
  \sup_{u\in\RR}\int_u^{u+1}
  \frac{|Z_{\scrP}(1+\varepsilon+it)|^2}
  {1+t^2}\dd t
  \ll_{\scrP,\Gamma}\frac1\varepsilon
  \qquad\text{for }0<\varepsilon\le\varepsilon_0.
\]
\end{enumerate}
\end{thm}

\begin{rmk}
If several $s_j$ are equal, the corresponding pole of $Z_{\scrP}(s)$ is still
at most simple.
\end{rmk}

The counting formulas below follow from Theorem~\ref{thm:main} by smoothing,
Mellin inversion, and a contour shift (cf., e.g.,
\cite[Sections~5.1 and~6.2]{MontgomeryVaughan}). Some of the details are given
in Appendix~\ref{appendix:contour}.

\begin{cor}
\label{cor:count}
For each $s_0\in\scrE$, let
\[
  c_{s_0}=\frac{1}{s_0}\Res_{s=s_0}Z_{\scrP}(s),
\]
with $c_{s_0}=0$ if $Z_{\scrP}$ is holomorphic at $s_0$. Then, as $T\to\infty$,
\begin{equation}
  \label{eq:count-cor}
  N_{\scrP}(T)=\sum_{s_0\in\scrE}c_{s_0}T^{s_0}
  +O\bigl(T^{(\deGa+1)/2}(\log T)^{1/4}\bigr).
\end{equation}
Moreover, for $X\ge2$,
\begin{equation}
  \label{eq:count-ms}
  \bigg(\frac1X\int_X^{2X}\Big|N_{\scrP}(T)-\sum_{s_0\in\scrE}c_{s_0}T^{s_0}\Big|^2\dd T\bigg)^{1/2}
  \ll X^{(\deGa+2)/3}(\log X)^{1/3}.
\end{equation}
\end{cor}
\begin{rmk}
The new feature of Corollary \ref{cor:count} compared with previous effective
circle counts is that it gives a Selberg-type counting formula, with the same
pointwise error exponent as in the counting theorem of Lax--Phillips for
$\Gamma$-orbits in $\HH^3$; see \cite{LaxPhillips}. In particular, the
exponent in the error term does not depend on the size of the spectral gap.
\end{rmk}
\begin{rmk}
\label{rmk:history}
\leavevmode
\begin{enumerate}[(1)]
\item Boyd \cite{Boyd73} proved that the packing exponent of a bounded
  Apollonian circle packing equals the Hausdorff dimension of its residual set,
  and subsequently showed \cite{Boyd82} that
\[
  N_{\scrP}(T)=T^{\deGa+o(1)},
  \qquad T\to\infty.
\]
Parker \cite{Parker} extended the relation between the circle-packing exponent
and the critical exponent to Kleinian circle packings. For the bounded packings
considered here, his result gives that $\deGa$ is the abscissa of convergence of
$Z_{\scrP}(s)$.

\item Kontorovich and Oh \cite{KontorovichOh} proved that for a bounded
  Apollonian circle packing, $N_{\scrP}(T)\sim c_{\scrP}T^{\deGa}$, using the
  equidistribution of expanding closed horospheres. Lee and Oh \cite{LeeOh}
  proved effective equidistribution of closed horospheres and obtained a power
  saving in the error term. Vinogradov \cite{Vinogradov} also proved an
  effective count using an estimate for the number of elements of $\Gamma$ in regions defined by the Cartan decomposition. In
  \cite[Paper IV, Section 5.6]{EdwardsThesis}, the dependence on the spectral
  gap in the horospherical approach of \cite{LeeOh} was improved.

\item Oh and Shah \cite{OhShah} established counting and equidistribution for
  more general families of circles invariant under Kleinian groups (see
  \cite[Section 2]{Oh} for a survey). Their setting also allows circles to
  intersect. Pan \cite{Pan} gave an effective version of the Oh--Shah results
  when $\deGa>1$ or $\Gamma$ is convex cocompact. The main input is exponential
  mixing of the geodesic flow on $\mathrm{T}^1(\scrM)$, proved in \cite{Stoyanov} for $\Gamma$
  convex cocompact, in \cite{MohammadiOh} when $\deGa>1$ (see also
  \cite{EdwardsOhDuke,Vinogradov}), and in \cite{LiPan} for all geometrically finite
  $\Gamma$ with cusps.
\item Kontorovich and Lutsko \cite{KontorovichLutsko} obtained an
  effective counting result for Kleinian sphere packings using the spectral
  theory of the Laplace operator. Kelmer, Kontorovich, and Lutsko
  \cite{KKL} proved that for the full symmetry group of an Apollonian
  packing, the base eigenvalue $\deGa(2-\deGa)$ is the only eigenvalue less than
  $1$. Their Corollary~1.1 gives the following effective count
\(  N_{\scrP}(T)=c_{\scrP}T^{\deGa}  +O\bigl(T^{(3\deGa+2)/5}(\log T)^{2/5}\bigr)\) in Apollonian
  circle packings. Corollary~\ref{cor:count} improves the exponent $(3\deGa+2)/5$ to
$(\deGa+1)/2$, and the power of $\log T$ from $2/5$ to $1/4$.

\item Pollicott \cite{Pollicott} gave a different proof of the asymptotic
  $N_{\scrP}(T)\sim c_{\scrP}T^{\deGa}$ for Apollonian circle packings by
  studying the Dirichlet series $Z_{\scrP}(s)$ and its meromorphic continuation
  using transfer operators. Pollicott and Urba\'nski \cite{PollicottUrbanski}
  developed this approach in a more general setting, obtaining counting
  asymptotics for conformal dynamical systems.
\end{enumerate}
\end{rmk}

\begin{rmk}
\label{rmk:simplifying-assumptions}
To simplify the notation, we assume from now on that $\Gamma$ is torsion-free
and that $\scrP=\Gamma C_0$ for some $C_0\in\scrP$. This gives
$\scrS=\Gamma_{C_0}\backslash\Pi_{C_0}$. Since $\infty\in\Omega(\Gamma)$, its
stabilizer $\Gamma_\infty$ is finite, and hence trivial under our assumption.

Our results hold without these assumptions: when $\scrP$ consists of several $\Gamma$-orbits, we integrate over each
component of $\scrS$ and add the integrals.
If $\Gamma$ has torsion or contains orientation-reversing isometries, we pass to a
torsion-free normal subgroup of finite index contained in $\PSL(2,\CC)$
and restrict to $\Gamma$-invariant functions.
On this cover we use the union of all packing-plane quotients, whose surface
measure is $\Gamma$-invariant, so only this invariant subspace contributes
to the residues. The group $\Gamma_\infty$ may then be non-trivial. We still sum over all of
$\Gamma$ in
\eqref{eq:eisenstein}, rather than over $\Gamma_{\infty}\backslash\Gamma$, so
\eqref{eq:period} still holds. See Remark~\ref{rmk:nontrivial-infty} for the
required change in Proposition~\ref{prop:zeta-pairing}.
\end{rmk}

\textbf{Organization:} In Section~\ref{Eperiod} we introduce normal coordinates
about a hyperbolic plane and prove the identity \eqref{eq:period}. In
Section~\ref{wave} we show that integration over $\scrS$ defines a bounded
functional on $H^1(\scrM)$, define its spectral measure $m_{\scrS}$, and prove
that $m_{\scrS}$ has bounded mass in intervals of length $2$ in the variable $\sqrt{\lambda-1}$
(Lemma~\ref{lem:spectral-window}). In Section~\ref{sec:continuation} we express
the integral of the Eisenstein series over $\scrS$ in terms of the resolvent of
$\scrL$, and prove Theorem~\ref{thm:main}\textit{(i)}.
Theorem~\ref{thm:main}\textit{(ii)} is proved in Section~\ref{sec:bounds}, and
Corollary~\ref{cor:count} in the appendix.\\

\textbf{Acknowledgements:} I would like to thank Michael Magee for explaining
the method of \cite{Pollicott} to me, and Hee Oh for useful discussions about
Eisenstein series in infinite volume during our joint work \cite{EO}, as well as
for helpful comments on an earlier draft.\\

\textbf{Declaration of LLM use:} ChatGPT (Astra GPT-6) was used to obtain the
wave equation argument in Sections \ref{wave} and \ref{sec:continuation},
replacing the author's longer resolvent kernel calculations. It also provided
the residue computation at $s=\deGa$ in the proof of Theorem
\ref{thm:main}\textit{(i)}. Claude (Opus 5.5 to Fable 5.2) produced Figure
\ref{fig:packings}, rewrote much of the author's first draft, and assisted in
checking the article for correctness and clarity.

\section{The Eisenstein period}
\label{Eperiod}

In this section we prove the identity \eqref{eq:period}. We first compute the
volume measure on $\HH^3$ in coordinates $(q,\rho)$, where $q$ is a point of a
hyperbolic plane $\Pi$, and $\rho$ is the signed distance from $\Pi$ along the
geodesic through $q$ normal to $\Pi$. After applying an isometry of $\HH^3$, we
may assume that $\Pi$ is the plane $\Pi=\{x+rj:x\in\RR,\ r>0\}$. For
$p=x+iy+rj$, let
\[
  h=\sqrt{y^2+r^2},\qquad
  y=h\tanh\rho,\qquad r=\frac{h}{\cosh\rho}.
\]
Then $q=x+hj$ is the closest point of $\Pi$ to $p$, and $\sinh\rho=y/r$. Since
\[
  d y^2+dr^2=d h^2+\frac{h^2}{\cosh^2\rho}\,d\rho^2,
\]
the hyperbolic metric is
\[
  \begin{aligned}
    \frac{d x^2+d y^2+d r^2}{r^2}
    =d\rho^2+\cosh^2\rho\frac{\,d x^2+d h^2}{h^2},
  \end{aligned}
\]
and $h^{-2}(dx^2+dh^2)$ is the hyperbolic metric on $\Pi$ in the coordinates
$x+hj$. Thus, for any hyperbolic plane $\Pi$,
\[
  \dd\vol(p)=\cosh^2\rho\,\dd\rho\,\dd S(q).
\]
For each $C\in\scrP$, we orient the geodesics normal to $\Pi_C$ towards $D_C$.
\begin{prop}
\label{prop:period}
The identity \eqref{eq:period} holds for $\Ree s\geq2$.
\end{prop}
\begin{proof}
Let $C$ have centre $c$ and radius $\mathsf{r}(C)$. In the upper half-space
model, we use spherical coordinates on $\Pi_C$: for $0\le\theta<\pi/2$,
$ 0\le\varphi<2\pi$, let
\[
  p=z+rj,\qquad
  z=c+\mathsf{r}(C)\sin\theta\,e^{i\varphi},\qquad r=\mathsf{r}(C)\cos\theta.
\]
In these coordinates, the restriction to $\Pi_C$ of the hyperbolic metric
$(|dz|^2+dr^2)/r^2$ is $(d\theta^2+\sin^2\theta\,d\varphi^2)/\cos^2\theta$, so
$dS=(\sin\theta/\cos^2\theta)\,\dd\theta\,\dd\varphi$.

Thus, for $\Ree s>1$,
\begin{equation}
  \label{eq:plane-height}
  \begin{aligned}
    \int_{\Pi_C}\htop(p)^s\dd S(p)
    &=\mathsf{r}(C)^s\int_0^{2\pi}\int_0^{\pi/2}
    (\cos\theta)^{s-2}\sin\theta\,\dd\theta\,\dd\varphi=\frac{2\pi}{s-1}\mathsf{r}(C)^s.
  \end{aligned}
\end{equation}
Since $\scrP=\Gamma C_0$, for $\Ree s\geq2$, we have
\begin{align*}
  \int_{\scrS}E(p,s)\dd S(p)
  &=\int_{\Gamma_{C_0}\backslash\Pi_{C_0}}
  \sum_{\gamma\in\Gamma}\htop(\gamma p)^s\dd S(p)\\
  &=\sum_{\gamma\in\Gamma/\Gamma_{C_0}}
  \int_{\Pi_{C_0}}\htop(\gamma p)^s\dd S(p)=\sum_{C\in\scrP}\int_{\Pi_C}\htop(p)^s\dd S(p)\\
  &=\frac{2\pi}{s-1}Z_{\scrP}(s)
\end{align*}
since the series and integrals above are absolutely convergent.
\end{proof}
\section{The spectral measure of the surface integral}
\label{wave}
In this section we show that integration over $\scrS$ is a bounded linear
functional on $H^1(\scrM)$, and we study its spectral measure $m_{\scrS}$ with
respect to $\scrL$. The main result is Lemma~\ref{lem:spectral-window}, which
bounds the mass of $m_{\scrS}$ in intervals of length $2$ in the variable $\sqrt{\lambda-1}$.

Let $H^1(\scrM)$ be the completion of $C_c^\infty(\scrM)$ with respect to the
norm
\[
  \norm{u}_{H^1(\scrM)}^2
  =\int_{\scrM}(|u|^2+|\nabla u|^2)\dd\vol,
\]
where $\nabla u$ is the hyperbolic gradient, and let $H^{-1}(\scrM)$ be its
dual. We consider the functional
\[
  \ell(u)=\int_{\scrS}u(q)\dd S(q),\qquad u\in C_c^\infty(\scrM).
\]

\begin{lem}
\label{lem:nearby-planes}
For every $p\in\HH^3$ and $a>0$,
\begin{equation}
  \label{eq:nearby-planes}
  \#\{C\in\scrP:\dist(p,\Pi_C)<a\}\le 1+e^{2a}.
\end{equation}
\end{lem}
\begin{proof}
Let $p$ be at the centre of the ball model. If $\dist(p,\Pi_C)<a$, then the disk
$D_C$ has spherical area at least $2\pi(1-\tanh a)$. Since the interiors of the
disks of $\scrP$ are disjoint,
\[
  \#\{C:\dist(p,\Pi_C)<a\}
  \le\frac{4\pi}{2\pi(1-\tanh a)}=1+e^{2a}.\qedhere
\]
\end{proof}

\begin{lem}
\label{lem:surface-trace}
For every $u\in C^1(\scrM)\cap H^1(\scrM)$,
\[
  \int_{\scrS}|u|^2\dd S\ll\norm{u}_{H^1(\scrM)}^2 .
\]
The extension $\ell\in H^{-1}(\scrM)$ satisfies
$\ell(u)=\int_{\scrS}u\dd S$ for every $u\in C^1(\scrM)\cap H^1(\scrM)$, where the integral converges absolutely.
\end{lem}

\begin{proof}
Let $C=C_0$, $\Pi=\Pi_C$, and $u\in C^1(\scrM)\cap H^1(\scrM)$.
We view $u$ as a $\Gamma$-invariant function on $\HH^3$.

Choose a fundamental domain $\mathcal D\subset\Pi$ for $\Gamma_C$. For $q\in\Pi$
and $\rho\in\RR$, let $q_{\rho}\in\HH^3$ denote the point at signed normal
distance $\rho$ from $\Pi$ whose closest point on $\Pi$ is $q$. Since $D_C$ is
the unique component of $\Omega(\Gamma)$ bounded by $C$, the group $\Gamma_C$
preserves $D_C$. We then have $(\gamma q)_\rho=\gamma(q_\rho)$ for all
$\gamma\in\Gamma_C$.

The one-dimensional Sobolev inequality gives, uniformly in $q$,
\[
  |u(q)|^2\ll\int_{-1}^1
  \bigl(|u(q_\rho)|^2+|\partial_\rho u(q_\rho)|^2\bigr)\dd\rho.
\]
Integrating over $\mathcal D$, using $|\partial_\rho u|\le |\nabla u|$ and the
volume measure from Section~\ref{Eperiod}, gives
\[
  \begin{aligned}
    \int_{\scrS}|u|^2\dd S
    &\ll\int_{\mathcal D}\int_{-1}^1
    \bigl(|u(q_\rho)|^2+|\nabla u(q_\rho)|^2\bigr)
    \cosh^2\rho\dd\rho\dd S(q)\\
    &=\int_{\Gamma_C\backslash\HH^3}
    \mathbf{1}_{\{\dist(p,\Pi)<1\}}(p)\cdot
    \bigl(|u(p)|^2+|\nabla u(p)|^2\bigr)\dd\vol(p).
  \end{aligned}
\]
The equality above holds because, in normal coordinates, $\mathcal D\times\RR$
is a fundamental domain for $\Gamma_C$ in $\HH^3$. Unfolding the last integral
from $\Gamma_C\backslash\HH^3$ to $\Gamma\backslash\HH^3$ gives
\[
  \int_{\scrM}\bigl(|u(p)|^2+|\nabla u(p)|^2\bigr)
  \sum_{\gamma\in\Gamma_C\backslash\Gamma}
  \mathbf{1}_{\{\dist(\gamma\tilde p,\Pi)<1\}}
  \dd\vol(p),
\]
where $\tilde p\in\HH^3$ is a lift of $p$. The sum counts the planes $\Pi_{C'}$,
$C'\in\scrP$, at distance less than $1$ from $\tilde p$, and is at most $1+e^2$
by Lemma~\ref{lem:nearby-planes}. This gives
\[
  \int_{\scrS}|u|^2\dd S\ll\norm{u}_{H^1(\scrM)}^2.
\]
Since $\operatorname{area}(\scrS)<\infty$, Cauchy--Schwarz gives
$|\ell(u)|\ll\norm{u}_{H^1(\scrM)}$ for $u\in C_c^\infty(\scrM)$. The extension to $H^1(\scrM)$
follows by
the density of $C_c^\infty(\scrM)$ in $H^1(\scrM)$.
\end{proof}

For an even real-valued $g\in C_c^\infty(\RR)$ and $p\in\HH^3$, define
\begin{equation}
  \label{eq:weighted-planes}
  W_g(p):=\sum_{C\in\scrP}\frac{g(\rho_C(p))}{\cosh\rho_C(p)},\qquad
  \rho_C(p)=\dist(p,\Pi_C).
\end{equation}
The sum is locally finite by Lemma~\ref{lem:nearby-planes} and is
$\Gamma$-invariant because $\scrP$ is. $W_g$ is therefore a well-defined function on $\scrM$.
If $g$ is supported in $(-a,a)$, then $W_g$ is supported in the
$a$-neighbourhood of $\scrS$. Lemma~\ref{lem:nearby-planes} gives
$|W_g|\le(1+e^{2a})\norm{g}_\infty$. The support of $W_g$ has volume at most $\operatorname{area}(\scrS)\int_{-a}^a\cosh^2\rho\dd\rho$, so
$W_g\in L^2(\scrM)$.

Laaksonen \cite[Section~3]{Laaksonen} uses similar sums on compact hyperbolic
$3$-manifolds and expands them in Laplace eigenfunctions. Since $\scrL$ on
$\scrM$ also has continuous spectrum, we use a point-pair kernel and the
spectral theorem instead.

For $g$ as above, we define
\begin{equation*}
  k_g(\varrho)=-\frac{g'(\varrho)}{2\pi\sinh \varrho},\qquad \varrho\in\RR.
\end{equation*}
The function $k_g$ is smooth at $\varrho=0$ and has compact support.
\begin{lem}
\label{lem:point-plane}
Let $\Pi\subset\HH^3$ be a hyperbolic plane, let $p\in\HH^3$, and let
$\rho=\dist(p,\Pi)$. Then
\begin{equation}
  \label{eq:point-plane}
  \int_\Pi k_g(\dist(p,q))\dd S(q)=\frac{g(\rho)}{\cosh\rho},
\end{equation}
and
\begin{equation}
  \label{eq:weighted-kernel}
  W_g(p)=\int_{\scrS}\left(\sum_{\gamma\in\Gamma}
  k_g\bigl(\dist(\tilde p,\gamma\tilde q)\bigr)\right)\dd S(q),
\end{equation}
where $\tilde p\in\HH^3$ is a lift of $p\in\scrM$, and $\tilde q\in\HH^3$ is a
lift of $q\in\scrS$.
\end{lem}
\begin{proof}
Let $q_0\in\Pi$ be the point closest to $p$, and let $v=\dist_\Pi(q_0,q)$. Then
\[
  \cosh\dist(p,q)=\cosh\rho\cosh v,
\]
so that $\sinh\varrho\,\dd\varrho=\cosh\rho\sinh v\,\dd v$, where
$\varrho=\dist(p,q)$. Polar coordinates on $\Pi$ then give
\begin{align*}
  \int_\Pi k_g(\dist(p,q))\dd S(q)
  &=2\pi\int_0^\infty k_g(\dist(p,q))\sinh v\,\dd v\\
  &=\frac{2\pi}{\cosh\rho}\int_\rho^\infty
  k_g(\varrho)\sinh\varrho\,\dd\varrho
  =\frac{g(\rho)}{\cosh\rho}.
\end{align*}
Writing $\scrS=\Gamma_{C_0}\backslash\Pi_{C_0}$ and unfolding the integral over
$\scrS$ on the right side of \eqref{eq:weighted-kernel} gives a sum of integrals
over $\gamma\Pi_{C_0}$, where $\gamma\in\Gamma/\Gamma_{C_0}$. Equation
\eqref{eq:point-plane} then gives the claimed formula for $W_g$.
\end{proof}

We next show that the kernel in \eqref{eq:weighted-kernel} is the kernel of a
function of $\scrL$. Since $\scrM$ is complete, the operator $\scrL$ with domain
$C_c^\infty(\scrM)$ is essentially self-adjoint on $L^2(\scrM)$
\cite[Theorem~2.4]{Str}; we denote its closure again by $\scrL$. The same
applies to $\HH^3$. For $\lambda\in\CC$, define
\[
  h_g(\lambda):=\int_{\RR}g(\tau)\cos(\tau\sqrt{\lambda-1})\dd\tau
  =-2\int_0^\infty g'(\tau)\,\frac{\sin(\tau\sqrt{\lambda-1})}{\sqrt{\lambda-1}}\dd\tau ,
\]
where the second expression is obtained by integration by parts. The functions involving
$\sqrt{\lambda-1}$ are defined by their power
series and are entire in $\lambda$.
The function $h_g$ is therefore entire. Since $g$ is smooth and compactly supported, $h_g$ is
real and bounded on $[0,\infty)$, and rapidly decreasing as $\lambda\to\infty$.
The spectral theorem therefore defines a bounded self-adjoint operator
$h_g(\scrL)$, and both formulas hold with $\scrL$ in place of $\lambda$
\cite[Section~VIII.3]{ReedSimon}.

For $u\in C_c^\infty(\scrM)$, the function
$U(\tau)=\frac{\sin(\tau\sqrt{\scrL-1})}{\sqrt{\scrL-1}}u$
is the unique solution in $C^2(\RR;L^2(\scrM))$, with values in
$\mathrm{Dom}(\scrL)$, of
\[
  U''+(\scrL-1)U=0,\qquad U(0)=0,\qquad U'(0)=u.
\]
For $u\in C_c^\infty(\HH^3)$, the corresponding solution on $\HH^3$ is
given by \cite[(2.3)]{MetcalfeTaylor}:
\[
  \Bigl(\frac{\sin(\tau\sqrt{\scrL-1})}{\sqrt{\scrL-1}}u\Bigr)(p)
  =\frac{1}{4\pi\sinh\tau}
  \int_{\{q\in\HH^3:\dist(p,q)=\tau\}}u(q)\dd A(q),\qquad \tau>0,
\]
where $dA$ is hyperbolic area measure.
Substituting this into the second formula for $h_g$ and using polar
coordinates gives the kernel $k_g(\dist(p,q))$ on $\HH^3$.
Uniqueness and finite propagation for the wave equation then give the kernel
on $\scrM$ by summing over $\Gamma$; see
\cite[Section~1, (1.1)--(1.4)]{CheegerGromovTaylor} and, for the constant shift,
\cite[Proposition~2.1]{MetcalfeTaylor}.
Thus, for $u\in C_c^\infty(\scrM)$ and $p\in\scrM$,
\begin{equation}
  \label{eq:automorphic-kernel}
  (h_g(\scrL)u)(p)
  =\int_{\scrM}\left(\sum_{\gamma\in\Gamma}
  k_g(\dist(\tilde p,\gamma\tilde q))\right)
  u(q)\,\dd\vol(q),
\end{equation}
where $\tilde p,\tilde q\in\HH^3$ are lifts of $p,q\in\scrM$, respectively.

We write $\langle\cdot,\cdot\rangle$ for the inner product on $L^2(\scrM)$. By
\cite[Corollary~2.6]{Str}, we have $\mathrm{Dom}(\scrL^{1/2})=H^1(\scrM)$ and
$\norm{u}_{H^1(\scrM)}^2=\norm{u}_2^2+\norm{\scrL^{1/2}u}_2^2$; in particular,
$\mathrm{Dom}(\scrL)\subset H^1(\scrM)$. Since the spectrum of $\scrL$ is
contained in $[\deGa(2-\deGa),\infty)$ and $\deGa(2-\deGa)>0$, $\scrL^{-1/2}$ is a bounded
operator from $L^2(\scrM)$ onto $H^1(\scrM)$.

By Lemma~\ref{lem:surface-trace} and the Riesz representation theorem,
there is a unique $f_0\in L^2(\scrM)$ such that
$\ell(\scrL^{-1/2}v)=\langle v,f_0\rangle$ for all $v\in L^2(\scrM)$, that is,
\begin{equation}
  \label{f0ident}
  \ell(u)=\langle \scrL^{1/2}u,f_0\rangle\qquad\text{for all }u\in H^1(\scrM).
\end{equation}
Let $E_\lambda$ denote the spectral resolution of $\scrL$
(cf.~\cite[Section~VIII.3]{ReedSimon}), and let $m_{\scrS}$ be the measure on
$\RR_{\geq0}$ given by
\[
  dm_{\scrS}(\lambda)=\lambda\dd\langle E_{\lambda}f_0,f_0\rangle .
\]
It is locally finite, since $m_{\scrS}([0,A])\le A\norm{f_0}_2^2$ for all
$A\ge0$.

Finally, we express $W_g$ in terms of $f_0$. Let $u\in C_c^\infty(\scrM)$. Then
$h_g(\scrL)u$ is in $C_c^\infty(\scrM)$ by \eqref{eq:automorphic-kernel}. Since
$W_g$ is real-valued and the kernel in \eqref{eq:automorphic-kernel} is real and
symmetric, \eqref{eq:weighted-kernel}, \eqref{eq:automorphic-kernel} and
\eqref{f0ident} give
\[
  \langle u,W_g\rangle=\ell\big(h_g(\scrL)u\big)
  =\big\langle\scrL^{1/2}h_g(\scrL)u,f_0\big\rangle
  =\big\langle u,\scrL^{1/2}h_g(\scrL)f_0\big\rangle .
\]
Here $\scrL^{1/2}h_g(\scrL)$ is bounded and self-adjoint because $h_g$ is real and rapidly decreasing on $[0,\infty)$. We now obtain
\begin{equation}
  \label{Wgident}
  W_g=\scrL^{1/2}h_g(\scrL)f_0
\end{equation}
as functions in $L^2(\scrM)$. By the decay of $h_g$, this function is in
$\mathrm{Dom}(\scrL)\subset H^1(\scrM)$.

\begin{lem}
\label{lem:pretrace}
For every even, real-valued $g\in C_c^\infty(\RR)$,
\[
  \int_0^{\infty}h_g(\lambda)\dd m_{\scrS}(\lambda)
  =\int_{\scrS}W_g(p)\dd S(p).
\]
\end{lem}
\begin{proof}
The function $W_g$ is smooth and, by \eqref{Wgident}, is in $H^1(\scrM)$.
Lemma~\ref{lem:surface-trace}, \eqref{f0ident} and the spectral theorem give
\[
  \int_{\scrS}W_g\dd S=\ell(W_g)
  =\langle\scrL h_g(\scrL)f_0,f_0\rangle
  =\int_0^\infty h_g(\lambda)\dd m_{\scrS}(\lambda).\qedhere
\]
\end{proof}

If $\supp g\subset(-a,a)$, then Lemma~\ref{lem:pretrace}, together with
\eqref{eq:weighted-planes} and \eqref{eq:nearby-planes}, gives
\begin{equation}
  \label{eq:trace-bound}
  \left|\int_{0}^{\infty} h_g(\lambda)\dd m_{\scrS}(\lambda)\right|
  \le\operatorname{area}(\scrS)(1+e^{2a})\norm{g}_\infty.
\end{equation}

For cocompact groups, the analogue of
$m_{\scrS}$ is the sum of the squared periods of the eigenfunctions over the
surface. For $T\ge1$, Laaksonen \cite[Lemma~3.3]{Laaksonen} proves that the sum
over eigenvalues $\lambda<1+T^2$ is $O(T)$. The corresponding bound on for bounded intervals in $\sqrt{\lambda}$, for closed
hypersurfaces in compact Riemannian manifolds, follows from the Kuznecov sum formula of Zelditch \cite{Zelditch}.

\begin{lem}
\label{lem:spectral-window}
For every $U\ge0$,
\begin{equation}
  \label{eq:spectral-window}
  m_{\scrS}\bigl(\{\lambda\in[1,\infty):|\sqrt{\lambda-1}-U|\le1\}\bigr)\ll1.
\end{equation}
\end{lem}
\begin{proof}
Choose a real, even function $b\in C_c^\infty((-1/2,1/2))$ whose Fourier
transform $\widehat{b}$, defined by
$\widehat b(\tau)=\int_{\RR}b(x)e^{-ix\tau}\dd x$, is non-negative on $\RR$ and
at least one on $[-1,1]$. Such a function can be obtained by convolving a
non-negative even smooth function supported in $(-1/4,1/4)$ with itself and
multiplying by a positive constant. For $U\ge0$, consider
\[
  g=(2b\cos(U\cdot))*(2b\cos(U\cdot)).
\]
Then $g$ is even, supported in $(-1,1)$, and $\norm{g}_\infty\ll1$ independently
of $U$. Moreover,
\[
  h_g(1+\tau^2)
  =\bigl(\widehat b(\tau-U)+\widehat b(\tau+U)\bigr)^2
  \ge1\qquad\text{for }|\tau-U|\le1.
\]
For $0\le\lambda<1$, we have
$h_g(\lambda)=\bigl(\int_{\RR}2b(x)\cos(Ux)\cosh(x\sqrt{1-\lambda})\dd x\bigr)^2\ge0$.
Therefore
\[
  \begin{aligned}
    m_{\scrS}\bigl(\{\lambda\in[1,\infty):|\sqrt{\lambda-1}-U|\le1\}\bigr)
    &\le\int_{0}^{\infty}h_g(\lambda)\dd m_{\scrS}(\lambda)\\
    &\ll\norm{g}_\infty\ll1,
  \end{aligned}
\]
by Lemma~\ref{lem:pretrace} and \eqref{eq:trace-bound}.
\end{proof}

\section{Continuation of the Dirichlet series}
\label{sec:continuation}

We now truncate the Eisenstein series near $\infty$ to obtain an $L^2$
function. The spectral theorem then gives the continuation of its integral over
$\scrS$.

For $p=z+rj\in\HH^3$, let $|p|=\sqrt{|z|^2+r^2}$ denote its Euclidean norm.
Since $\infty\in\Omega(\Gamma)$ and $\Gamma_\infty$ is trivial, we may choose
$R>0$ such that every packing plane is contained in the set
$\{p\in\HH^3\, : \, |p|<R/2\}$ and $ \{p\in\HH^3\, : \, |p|>R\}$ projects
injectively into $\scrM$.

We fix a smooth non-decreasing function $\chi_0:[0,\infty)\to\RR$ such that
\[
  \chi_0(t)=
  \begin{cases}
    0\quad &\text{ for }t\le R\\
    1 \quad &\text{ for }t\ge 2R.
  \end{cases}
\]

Define $\chi:\HH^3\to\RR$ by $\chi(p)=\chi_0(|p|)$, and let
\begin{equation}
  \label{eq:eta-source}
  \eta_s=(\scrL-\lambda_s)(\chi\,\htop^s).
\end{equation}
Since $(\scrL-\lambda_s)\htop^s=0$, a direct calculation gives
\[
  \eta_s(p)=r^{s+2}\left(-\Delta_{\euc}\chi(p)
  +(1-2s)\frac{\chi_0'(|p|)}{|p|}\right),\qquad p=z+rj\in\HH^3.
\]
The function $\eta_s$ vanishes outside $R<|p|<2R$, which projects injectively
into $\scrM$. We extend it by zero to $\scrM$. On this region,
\[
  |\eta_s(p)|\ll(1+|s|)r^{\Ree s+2}.
\]
The same bound, with powers of $|\log r|$, applies to its derivatives with respect to $s$. Since
$d\vol=r^{-3}\,dx\,dy\,dr$, it follows that $s\mapsto\eta_s$ is
holomorphic with values in $L^2(\scrM)$ for $\Ree s>0$, and
\begin{equation}
  \label{eq:eta-bound}
  \|\eta_s\|_2\ll1+|s|\qquad(1\le\Ree s\le2).
\end{equation}

\begin{lem}
\label{lem:Ecuttoff}
For $\Ree s>\deGa$, let $u_s:=E(\cdot,s)-\chi\htop^s$, where $\chi\htop^s$ is
viewed as a function on $\scrM$, extended by zero outside the image of
$\{p\in\HH^3\,:\,|p|>R\}$. Then $u_s$ is in
$\mathrm{Dom}(\scrL)\subset H^1(\scrM)$, and
\begin{equation}
  \label{usdef}
  u_s=-\big(\scrL-\lambda_s\big)^{-1}\eta_s.
\end{equation}
\end{lem}
\begin{proof}
Let $\sigma=\Ree s>\deGa$. The spectral theorem gives
\[
  (\scrL-\lambda_s)^{-1}
  =\int_0^\infty e^{-(s-1)\tau}
  \frac{\sin(\tau\sqrt{\scrL-1})}{\sqrt{\scrL-1}}\dd\tau.
\]
The integral converges in operator norm, since the norm of
$\sin(\tau\sqrt{\scrL-1})/\sqrt{\scrL-1}$ is at most
$\sinh((\deGa-1)\tau)/(\deGa-1)$ on $L^2(\scrM)$ and at most $\tau$ on
$L^2(\HH^3)$. Integrating the wave formula from Section~\ref{wave} gives
the resolvent kernel
\[
  G_s(p,q)=\frac{e^{-(s-1)\dist(p,q)}}{4\pi\sinh\dist(p,q)}
\]
on $\HH^3$, and its sum over $\Gamma$ on $\scrM$. In particular, these kernels
are positive when $s=\sigma$ is real.

Let $v_s=(1-\chi)\htop^s$ on $\HH^3$. Since
\[
  \|v_s\|_2^2\ll
  \int_0^{2R}r^{2\sigma-3}\dd r<\infty,
  \qquad
  (\scrL-\lambda_s)v_s=-\eta_s\in L^2(\HH^3),
\]
we have $v_s\in\mathrm{Dom}(\scrL_{\HH^3})$ by
\cite[Theorem~2.4 and Corollary~2.5]{Str}. This gives
\[
  v_s(p)=-\int_{\HH^3}G_s(p,q)\eta_s(q)\dd\vol(q).
\]
Since $|G_s|=G_\sigma$ and $\supp\eta_s$ projects injectively to $\scrM$, we
have, almost everywhere on $\scrM$,
\[
  \sum_{\gamma\in\Gamma}|v_s(\gamma\tilde p)|
  \le\big((\scrL-\lambda_\sigma)^{-1}|\eta_s|\big)(p),
\]
where $\tilde p$ is a lift of $p$. The right-hand side is in $L^2(\scrM)$, so the sum on the left defines an $L^2$ function. Moreover,
$\sum_{\gamma\in\Gamma}v_s(\gamma\tilde p)=u_s(p)$, since
$\sum_{\gamma\in\Gamma}(\chi\htop^s)(\gamma\tilde p)$ is the extension of
$\chi\htop^s$ from the exterior region to $\scrM$, and so $u_s\in L^2(\scrM)$.

Finally, $u_s\in L^2(\scrM)$ and
$(\scrL-\lambda_s)u_s=-\eta_s\in L^2(\scrM)$ as distributions. This implies that $u_s\in\mathrm{Dom}(\scrL)$ by \cite[Theorem~2.4 and Corollary~2.5]{Str}, and \eqref{usdef} follows.
\end{proof}

\begin{rmk}
\label{rmk:nontrivial-infty}
If $|\Gamma_\infty|>1$, after applying a translation, we may assume that
$\Gamma_\infty$ acts by rotations and reflections fixing $z=0$.
Then $\chi$ is $\Gamma_\infty$-invariant, the image of $\{|p|>R\}$ in
$\scrM$ is its quotient by $\Gamma_\infty$, and $\eta_s$ is defined on
this quotient. Lemma~\ref{lem:Ecuttoff} holds with
\[
  u_s=E(\cdot,s)-|\Gamma_{\infty}|\chi\htop^s
  =-|\Gamma_{\infty}|(\scrL-\lambda_s)^{-1}\eta_s.
\]
Thus the final expression in \eqref{eq:zeta-pairing} and the formulas below for $Z_{\scrP}$ and its residues gain a factor $|\Gamma_\infty|$.
The right-hand side of \eqref{eq:Zbound} acquires the factor
$|\Gamma_\infty|^2$.
\end{rmk}

\begin{prop}
\label{prop:zeta-pairing}
For $\Ree s>\deGa$,
\begin{equation}
  \label{eq:zeta-pairing}
  \int_{\scrS}E(p,s)\,dS(p)=\ell(u_s) =-\ell\big((\scrL-\lambda_s)^{-1}\eta_s\big).
\end{equation}
\end{prop}
\begin{proof}
The function $\chi\htop^s$ on $\scrM$ vanishes on $\scrS$, since every packing
plane is contained in $\{p\in\HH^3\,:\,|p|<R/2\}$. By Lemma~\ref{lem:Ecuttoff},
$u_s$ is a smooth function in $\mathrm{Dom}(\scrL)\subset H^1(\scrM)$. 
Lemma~\ref{lem:surface-trace} and \eqref{usdef} now give
\begin{align*}
  \int_{\scrS}E(p,s)\dd S(p)
  &=\int_{\scrS}\big(E(p,s)-\chi(p)\htop(p)^s\big)\dd S(p)\\
  &=\ell(u_s)=-\ell\big((\scrL-\lambda_s)^{-1}\eta_s\big),
\end{align*}
and the integrals are absolutely convergent.
\end{proof}

\begin{proof}[Proof of Theorem~\ref{thm:main}\textit{(i)}]
For $\Ree s>1$ and $s\notin\scrE$, the spectral theorem gives
\[
  (\scrL-\lambda_s)^{-1}
  =\int_{\spec(\scrL)}\frac{1}{\lambda-\lambda_s}\dd E_\lambda.
\]
The quotient $\sqrt{1+\lambda}/|\lambda-\lambda_s|$ is bounded
uniformly for $\lambda\in\spec(\scrL)$ and for $s$ in any compact
subset of the half-plane $\Ree s>1$ disjoint from $\scrE$. The resolvent is therefore holomorphic there as an operator from $L^2(\scrM)$ to $H^1(\scrM)$, and has simple poles at the points of $\scrE$; see \cite[Section~VIII.3]{ReedSimon}.

Since $\eta_s$ is holomorphic with values in $L^2(\scrM)$ and $\ell$ is bounded
on $H^1(\scrM)$, Propositions~\ref{prop:period} and~\ref{prop:zeta-pairing} give
the required continuation of $Z_{\scrP}$.

We next compute the residues. Since $\lambda_{s_0}-\lambda_s=(s-s_0)(s+s_0-2)$,
the residue of $(\scrL-\lambda_s)^{-1}$ at $s=s_0$ is
\[
  \frac{1}{2(s_0-1)}P_{\lambda_{s_0}},
\]
where $P_{\lambda_{s_0}}$ is the orthogonal projection onto the
$\lambda_{s_0}$-eigenspace of $\scrL$. Thus
\[
  \Res_{s=s_0} Z_{\scrP}(s)=\frac{-1}{4\pi}\ell\big(P_{\lambda_{s_0}}\eta_{s_0}\big).
\]

If $\lambda_{s_0}$ is a simple eigenvalue with real unit eigenfunction $\psi$,
Lemma~\ref{lem:surface-trace} gives
\begin{equation}
  \label{eq:residue}
  \Res_{s=s_0}Z_{\scrP}(s)
  =-\frac1{4\pi}
  \left(\int_{\scrS}\psi\dd S\right)
  \left(\int_{\scrM}\eta_{s_0}\psi\dd\vol\right).
\end{equation}
For an eigenvalue of higher multiplicity, the right-hand side is summed over a
real orthonormal basis of the eigenspace. This residue may be zero.

It remains to show that the residue at $s=\deGa$ is positive. For the base
eigenvalue $\lambda_{\deGa}$, let $\phi_0$ be the positive eigenfunction
normalized by $\norm{\phi_0}_2=1$.

We show that $\int_{\scrM}\eta_{\deGa}\phi_0\dd\vol<0$. On
$\{p\in\HH^3:|p|>R\}$, which we identify with its image in $\scrM$, define
\[
  w(p)=\frac{\htop(p)^{\deGa}}{|p|^{2\deGa}}
  \int_R^{|p|}t^{2\deGa}\chi_0'(t)\dd t,
\]
and extend $w$ by zero to the rest of $\scrM$. The function $w$ is smooth.
Differentiating $w$ and using the calculation following
\eqref{eq:eta-source} gives
\[
  \bigl(\scrL-\lambda_{\deGa}\bigr)w(p)
  =\eta_{\deGa}(p)
  +2\deGa\frac{\htop(p)^{\deGa+2}}{|p|}\chi_0'(|p|),
\]
where the last term is also extended by zero to $\scrM$.

Since $\chi_0'(t)=0$ for $t\ge2R$, the integral in the definition of $w$ is
bounded independently of $p$, and so $|w(p)|\ll(\htop(p)/|p|^2)^{\deGa}$ for
$|p|>R$, and
\[
  \begin{aligned}
    &\|w\|_2^2\ll\int_R^\infty t^{-2\deGa-1}\dd t
    \int_0^{\pi/2}(\cos\theta)^{2\deGa-3}\sin\theta\dd\theta
    <\infty,
  \end{aligned}
\]
since $\deGa>1$. The additional term in the formula for
$(\scrL-\lambda_{\deGa})w$ is in $L^2(\scrM)$, since it is
$O(\htop(p)^{\deGa+2})$ and vanishes outside $R<|p|<2R$. Thus
$w\in\mathrm{Dom}(\scrL)$ by \cite[Theorem~2.4 and Corollary~2.5]{Str}. Taking
the inner product with $\phi_0$ and using $\scrL\phi_0=\lambda_{\deGa}\phi_0$
gives
\[
  \begin{aligned}
    \int_{\scrM}\eta_{\deGa}\phi_0\dd\vol
    &=-2\deGa\int_{\{p\in\HH^3:R<|p|<2R\}}
    \frac{\htop(p)^{\deGa+2}}{|p|}
    \chi_0'(|p|)\phi_0(p)\dd\vol(p)<0.
  \end{aligned}
\]
Here $\phi_0$ denotes its lift to $\HH^3$ in the integral on the right. The
strict inequality follows from $\phi_0>0$ and $\chi_0'\ge0$, with $\chi_0'$ not
identically zero. Applying \eqref{eq:residue} to the base eigenvalue now gives
\begin{equation}
  \label{eq:leading-residue}
  \Res_{s=\deGa}Z_{\scrP}(s)
  =-\frac1{4\pi}
  \left(\int_{\scrS}\phi_0\dd S\right)
  \left(\int_{\scrM}\eta_{\deGa}\phi_0\dd\vol\right)>0.
\end{equation}

This proves part \textit{(i)} of Theorem~\ref{thm:main}.
\end{proof}
\begin{rmk}
Alternatively, Parker's identification of the abscissa of convergence
\cite{Parker}, together with Landau's theorem
\cite[Lemma~15.1]{MontgomeryVaughan}, implies that $\deGa$ is a pole. Its
residue is then positive because $Z_{\scrP}(s)>0$ for real $s>\deGa$.
\end{rmk}

\section{Bounds on \texorpdfstring{$Z_{\scrP}(s)$}{the Dirichlet series}}
\label{sec:bounds}

For $\Ree s>1$ and $s\notin\scrE$, let
\begin{equation}
  \label{eq:hs}
  h_s(\lambda)=\frac1{\lambda-\lambda_s},\qquad
  \lambda\in\spec(\scrL).
\end{equation}
By Propositions~\ref{prop:period} and~\ref{prop:zeta-pairing} and continuation,
\[
  Z_{\scrP}(s)=-\frac{s-1}{2\pi}\ell(h_s(\scrL)\eta_s).
\]
Equation~\eqref{f0ident} and the spectral theorem give
\[
  |\ell(h_s(\scrL)\eta_s)|^2
  =|\langle\eta_s,\scrL^{1/2}h_{\bar s}(\scrL)f_0\rangle|^2
  \le\|\eta_s\|_2^2\int_0^\infty|h_s(\lambda)|^2\dd m_{\scrS}(\lambda).
\]
We now have
\begin{equation}
  \label{eq:Zbound}
  |Z_{\scrP}(s)|^2\le\frac{\|\eta_s\|_2^2}{4\pi^2}
  \int_0^\infty|s-1|^2|h_s(\lambda)|^2\dd m_{\scrS}(\lambda).
\end{equation}

\begin{proof}[Proof of Theorem~\ref{thm:main}\textit{(ii)}]
Choose $0<\varepsilon_0<\min_{s_0\in\scrE}(s_0-1)$. Let $s=1+\varepsilon+it$,
where $0<\varepsilon\le1$, and suppose that either $\varepsilon\le\varepsilon_0$
or $|t|\ge1$.

For an eigenvalue $\lambda_{s_0}<1$,
\[
  |s-1|\,|h_s(\lambda_{s_0})|
  =\frac{|s-1|}{|s-s_0|\,|s+s_0-2|}
  \le\frac1{|s-s_0|}.
\]
In either case, $|s-s_0|\ge\min\{1,s_0-1-\varepsilon_0\}>0$. Since $\scrE$ is finite and
$m_{\scrS}$ is locally finite, the eigenvalues below $1$ contribute $O(1)$ to
the integral in \eqref{eq:Zbound}.

For $\lambda\ge1$, let $\tau=\sqrt{\lambda-1}$. Since
$\lambda-\lambda_s=\tau^2+(s-1)^2$, we have
\begin{equation}
  \label{eq:above1}
  |s-1|^2|h_s(\lambda)|^2
  =\frac{\varepsilon^2+t^2}
  {(\varepsilon^2+(t+\tau)^2)(\varepsilon^2+(t-\tau)^2)}
  \le\frac1{\varepsilon^2+(\tau-|t|)^2}.
\end{equation}
By Lemma~\ref{lem:spectral-window}, the set $|\tau-|t||\le1$ and each of the
sets $k<|\tau-|t||\le k+1$, $k\ge1$, have bounded $m_{\scrS}$-mass. This gives
\begin{equation}
  \label{eq:above1int}
  \int_1^\infty|s-1|^2|h_s(\lambda)|^2\dd m_{\scrS}(\lambda)
  \ll\varepsilon^{-2}+\sum_{k\ge1}k^{-2}
  \ll\varepsilon^{-2}.
\end{equation}
For $|t|\ge1$, \eqref{eq:eta-bound} and \eqref{eq:Zbound} now give
$|Z_{\scrP}(1+\varepsilon+it)|\ll |t|/\varepsilon$.

Now suppose $\varepsilon\le\varepsilon_0$. By \eqref{eq:eta-bound},
\eqref{eq:Zbound} and \eqref{eq:above1},
\[
  \frac{|Z_{\scrP}(1+\varepsilon+it)|^2}{1+t^2}
  \ll1+\int_1^\infty
  \frac{1}{\varepsilon^2+(\tau-|t|)^2}\dd m_{\scrS}(\lambda).
\]
Integrate over $u\le t\le u+1$. For fixed $\tau$,
\[
  \int_u^{u+1}\frac{\dd t}{\varepsilon^2+(\tau-|t|)^2}
  \le\frac{2\pi}{\varepsilon}.
\]
The contribution from $|\tau-|u||\le2$ is $O(\varepsilon^{-1})$, since this
set has bounded $m_{\scrS}$-mass. For $k<|\tau-|u||\le k+1$, $k\ge2$, the
integral in $t$ is at most $(k-1)^{-2}$, and the $m_{\scrS}$-mass is again
bounded. Thus, uniformly in $u$,
\[
  \int_u^{u+1}\frac{|Z_{\scrP}(1+\varepsilon+it)|^2}{1+t^2}\dd t
  \ll1+\varepsilon^{-1}+\sum_{k\ge2}(k-1)^{-2}
  \ll\varepsilon^{-1}.\qedhere
\]
\end{proof}

\appendix
\section{Contour integrals}
\label{appendix:contour}

We apply the argument of \cite[Appendix~C.2, pp.~448--452]{Edwards}, with $T$
and $Z_{\scrP}$ in place of $M$ and $L_k$. We describe the changes needed to
prove Corollary~\ref{cor:count}.

\begin{proof}[Proof of Corollary~\ref{cor:count}]
Let $\psi_U$ and $\Psi_U$ be as in \cite{Edwards}, with $U>2$ and
$0\le\psi_U\le1$. Thus $\psi_U(x)=1$ for $x\le1-U^{-1}$ and $\psi_U(x)=0$ for
$x\ge1+U^{-1}$. Let $\varepsilon_0$ be as in Theorem~\ref{thm:main}\textit{(ii)},
decreased if necessary so that $1+\varepsilon_0<\min\scrE$.
For $0<\varepsilon\le\varepsilon_0$, let
\[
  I(T)=\frac1{2\pi i}\int_{1+\varepsilon+i\RR}
  Z_{\scrP}(s)\Psi_U(s)T^s\dd s.
\]
The pointwise bound in Theorem~\ref{thm:main}\textit{(ii)}, together with
\cite[(81)]{Edwards} with $c=3$, allows us to move the contour from
$2+i\RR$ to $1+\varepsilon+i\RR$, so \cite[Lemma~37]{Edwards} is not needed. The poles crossed
by the contour are at most simple, and
$\Psi_U(\zeta)=\zeta^{-1}+O(U^{-1})$ for $\zeta\in\scrE$, so
\cite[(83)]{Edwards} becomes
\begin{equation}
  \label{eq:smoothcount}
  \sum_{C\in\scrP}\psi_U\Bigl(\frac1{T\mathsf r(C)}\Bigr)
  =\sum_{\zeta\in\scrE}c_\zeta T^\zeta+I(T)
  +O\Bigl(\frac{T^{\deGa}}U\Bigr).
\end{equation}

Let $s=1+\varepsilon+it$. The cases $c=1$ and $c=3$ of \cite[(81)]{Edwards} give
$\sqrt{1+t^2}|\Psi_U(s)|\ll\min\{1,U^2(1+|t|)^{-2}\}$. Applying the mean-square bound
in Theorem~\ref{thm:main}\textit{(ii)} on the intervals $[n,n+1]$ and summing gives
\begin{equation}
  \label{eq:Ibound}
  |I(T)|\ll T^{1+\varepsilon}U\varepsilon^{-1/2},
  \qquad
  \int_{\RR}|Z_{\scrP}(s)\Psi_U(s)|^2\dd t\ll U\varepsilon^{-1}.
\end{equation}
These bounds replace \cite[(84)]{Edwards}.
Their dependence on $\varepsilon$ gives the logarithmic factors below.

Let $F_U(T)$ denote the sum on the left-hand side of \eqref{eq:smoothcount}.
The inequalities
\[
  F_U\!\left(\frac{T}{1+U^{-1}}\right)
  \le N_{\scrP}(T)
  \le F_U\!\left(\frac{T}{1-U^{-1}}\right)
\]
give
\[
  \Bigl|N_{\scrP}(T)-\sum_{\zeta\in\scrE}c_\zeta T^\zeta\Bigr|
  \ll T^{1+\varepsilon}U\varepsilon^{-1/2}+T^{\deGa}/U.
\]
Taking $\varepsilon=(\log T)^{-1}$ and $U=T^{(\deGa-1)/2}(\log T)^{-1/4}$ proves
\eqref{eq:count-cor}.

For \eqref{eq:count-ms}, keep $\varepsilon$ and $U$ fixed on $X\le T\le2X$.
Plancherel's theorem, after the substitution $T=e^x$, and \eqref{eq:Ibound} give
\[
  \begin{aligned}
    \int_{X/2}^{4X}|I(T)|^2\dd T
    &\ll X^{3+2\varepsilon}\int_{\RR}
    |Z_{\scrP}(1+\varepsilon+it)\Psi_U(1+\varepsilon+it)|^2\dd t\\
    &\ll X^{3+2\varepsilon}U/\varepsilon.
  \end{aligned}
\]
The comparison at $T/(1\pm U^{-1})$ bounds the left-hand side of
\eqref{eq:count-ms} by
\[
  X^{1+\varepsilon}(U/\varepsilon)^{1/2}+X^{\deGa}/U.
\]
Taking $\varepsilon=(\log X)^{-1}$ and $U=X^{2(\deGa-1)/3}(\log X)^{-1/3}$
proves \eqref{eq:count-ms}. For sufficiently large $T$ or $X$, these choices satisfy $U>2$ and
$0<\varepsilon\le\varepsilon_0$. Enlarging the implied constants covers the
remaining bounded ranges.
\end{proof}

\end{document}